\documentclass[12pt]{article}
\usepackage{hyperref,amsfonts,amssymb}
\usepackage{amsfonts,amssymb}
\usepackage{mathtools}
\usepackage{rotating}

\def\C{\mathbb{C}}

\def\N{\mathbb{N}}
\def\Z{\mathbb{Z}}

\def\bq{ \begin{equation} }
\def\eq{ \end{equation} }
\def\ben{ \begin{eqnarray} }
\def\en{ \end{eqnarray} }

\def\frac#1#2{{#1\over #2}}

\def\on#1#2{\mathop{\vbox{\ialign{##\crcr\noalign{\kern2pt}
$\scriptstyle{#2}$\crcr\noalign{\kern2pt\nointerlineskip}
\kern-2pt$\hfil\displaystyle{#1}\hfil$\crcr}}}\limits}

\begin{document}

\title{Poisson structures on loop spaces of $\C P^n$ and an $r$-matrix associated with the universal elliptic curve}
\author{Alexander Odesskii}
   \date{}
\vspace{-20mm}
   \maketitle
\vspace{-7mm}
\begin{center}
Brock University, 1812 Sir Isaac Brock Way, St. Catharines, ON, L2S 3A1 Canada\\[1ex]
email: \\
\texttt{aodesski@brocku.ca}
\end{center}

\medskip

\begin{abstract}

\medskip

We construct a family of Poisson structures of hydrodynamic type on the loop space of $\C P^{n-1}$. This family is parametrized by the moduli space of elliptic curves or, in other words, by the modular parameter $\tau$. This family can be lifted to a homogeneous Poisson structure on the loop space of $\C^n$ but in order to do that we need to upgrade the modular parameter $\tau$ to an additional field $\tau(x)$ with Poisson brackets $\{\tau(x),\tau(y)\}=0,~~\{\tau(x),z_a(y)\}=2\pi i~ z_a(y)~\delta^{\prime}(x-y)$ where $z_1,...,z_n$ are coordinates on $\C^n$.  These homogeneous Poisson structures can be written in terms of an elliptic $r$-matrix of hydrodynamic type.

\medskip

\end{abstract}

\newpage
\tableofcontents
\newpage

\section{Introduction}

Consider a holomorphic Poisson structure on a complex manifold $M$. It can be written in local coordinates $z_1,...,z_n$ as 
$$\{z_a,z_b\}=g_{ab}$$
where $g_{ab}$ are holomorphic functions in $z_1,...,z_n$. Such Poisson structures  are well understood locally \cite{lie,dar}: in a neighborhood of any generic point one can always choose coordinates such that $g_{ab}$ 
are constants. It is an interesting and hard problem however to classify and/or study global holomorphic Poisson structures \cite{bond, pol, pym} on a given compact complex manifold $M$. Even in a basic case $M=\C P^n$ no classification results are known for $n>3$. 

Consider the simplest case of $\C P^2$. Let $p_1,p_2$ be affine coordinates on $\C P^2$. It is clear that 
$$\{p_1,p_2\}=Q(p_1,p_2)$$ 
where $Q$ is a polynomial in 
$p_1,p_2$. Moreover, $\{\frac{1}{p_1},\frac{p_2}{p_1}\}=-\frac{Q(p_1,p_2)}{p_1^3}$ \Big(resp. $\{\frac{1}{p_2},\frac{p_1}{p_2}\}=\frac{Q(p_1,p_2)}{p_2^3}$\Big) should be a polynomial in $\frac{1}{p_1},\frac{p_2}{p_1}$ \Big(reps. in $\frac{1}{p_2},\frac{p_1}{p_2}$\Big). Therefore, $Q$ is an arbitrary polynomial of degree $3$. After a projective change of variables we obtain (in generic case)
\begin{equation}\label{pr1}
\{p_1,p_2\}=p_1^2-4p_2^3+g_2p_2+g_3 
\end{equation}
where $g_2,~g_3$ are parameters. In other words, Poisson structures on $\C P^2$ are parametrized by it's set of zeros which is an elliptic curve (or it's degeneration) given by a cubic in $\C P^2$.

Let $z_1:z_2:z_3$ be homogeneous coordinates on $\C P^2$. It turns out that the Poisson structure (\ref{pr1}) on $\C P^2$ can be lifted to a homogeneous Poisson structure on $\C^3$ as follows:
\begin{equation}\label{hm1}
\{z_1,z_2\}=\frac{\partial \tilde{Q}}{\partial z_3},~~~ \{z_2,z_3\}=\frac{\partial \tilde{Q}}{\partial z_1},~~~\{z_3,z_1\}=\frac{\partial \tilde{Q}}{\partial z_2}
\end{equation} 
where $\tilde{Q}=\frac{1}{3} (z_1^2z_3-4z_2^3+g_2z_2z_3^2+g_3 z_3^3)$. Indeed, computing Poisson bracket between $p_1=\frac{z_1}{z_3}$ and $p_2=\frac{z_2}{z_3}$ by virtue of (\ref{hm1}) we obtain (\ref{pr1}).

In general, given a homogeneous Poisson structure on $\C^n$ of the form
\begin{equation}\label{pg}
\{z_a,z_b\}=\sum_{1\leq c,d\leq n} q_{ab}^{cd}~ z_cz_d 
\end{equation}
one can always obtain a Poisson structure on $\C P^{n-1}$ with homogeneous coordinates $z_1:...:z_n$. Moreover, it is known \cite{bond, pol} that any Poisson structure on $\C P^{n-1}$ can be obtained in this way. 

In this paper we are interested in Poisson structures of the form
\begin{equation}\label{phg}
\{z_a(x),z_b(y)\}=g_{ab}~\delta^{\prime}(x-y)+\sum_{1\leq c\leq n}h^c_{ab}z_c^{\prime}(x)~\delta(x-y)
\end{equation}
where $z_1,...,z_n$ are local coordinates on a complex manifold\footnote{These Poisson structures appeared in the literature by different names such as coisson structures \cite{bdr}, Poisson structures on loop space of $M$ \cite{moh}, Poisson structures of hydrodynamic type \cite{dn}. In this paper we often abuse language by referring to these as to just Poisson structures on $M$.} $M$ promoted to fields $z_1(x),...,z_n(x)$ and $g_{ab},~h^c_{ab}$ are holomorphic functions in $z_1(x),...,z_n(x)$. These Poisson structures are also well understood locally in non-degenerate case \cite{dn}: if $\det(g_{ab})\ne 0$, then one can always choose local coordinates such that $g_{ab}$ are constants and $h^c_{ab}=0$. In this paper however we are interested in global holomorphic Poisson structures of the form (\ref{phg}) on projective spaces $\C P^n$.

Consider the simplest case of $\C P^1$ with affine coordinate $p$. We have 
$$\{p(x),p(y)\}=G(p(x))\delta^{\prime}(x-y)+H(p(x))p^{\prime}(x)\delta(x-y)$$ 
where $G(p),~H(p)$ are polynomials. Anticommutativity and Jacobi identity lead to the constraint $H=\frac{1}{2}G^{\prime}$. Computing $\{\frac{1}{p(x)},\frac{1}{p(y)}\}$ we obtain another constraint: $\frac{G(p)}{p^4}$ should be a polynomial in $\frac{1}{p}$. 
Therefore, $G(p)$ is an arbitrary polynomial of degree less or equal to 4 and our Poisson structure reads:
\begin{equation}\label{phpr1}
\{p(x),p(y)\}=G(p(x))\delta^{\prime}(x-y)+\frac{1}{2}G^{\prime}(p(x))p^{\prime}(x)\delta(x-y). 
\end{equation}
Let $z_1:z_2$ be homogeneous coordinates on $\C P^1$. Is it possible to lift the Poisson structure (\ref{phpr1}) on $\C P^1$ to a homogeneous Poisson structure on $\C^2$? Such Poisson structure  should have  the form
\begin{equation}\label{hom}
\{z_a(x),z_b(y)\}=\sum_{1\leq c,d\leq 2}q_{ab}^{cd}z_c(x)z_d(x)\delta^{\prime}(x-y)+\sum_{1\leq c,d\leq 2}r_{ab}^{cd}z_c(x)z^{\prime}_d(x)\delta(x-y) 
\end{equation}
where $q_{ab}^{cd},~r_{ab}^{cd}$ are constants. The answer to this question is negative:

{\bf Proposition 1.} If a polynomial $G(p)$ has distinct roots, then there is no Poisson structure of the form (\ref{hom}) such that $p(x)=\frac{z_1(x)}{z_2(x)}$ satisfy (\ref{phpr1}) by virtue of (\ref{hom}).

{\bf Proof.} Let roots of $G(p)$ be distinct. Using projective transformations we can send these roots to $0,~1,~s,~\infty$. Therefore, without loss of generality, we can set $G(p)=p(p-1)(p-s)$ where 
$s\ne 0,~1$. Computing $\{p(x),p(y)\}=\{\frac{z_1(x)}{z_2(x)},\frac{z_1(y)}{z_2(y)}\}$ by virtue of (\ref{hom}) and comparing the result with (\ref{phpr1}) we obtain linear equations for coefficients 
$q_{ab}^{cd},~r_{ab}^{cd}$. Anticommutativity and Jacobi identity for (\ref{hom}) give additional linear and quadratic equations for $q_{ab}^{cd},~r_{ab}^{cd}$. Direct calculation shows that this system 
of equations does not have any solutions if $s\ne 0,~1$. $\Box$
 
There exists however a way to construct a homogeneous Poisson structure similar to (\ref{hom}) but in order to do that we need to introduce a dynamical variable. This can be done in a natural way: the family (\ref{phpr1}) is parametrized by an isomorphism class of elliptic curve\footnote{Which is a double cover of $\C P^1$ without four roots of $G(p)$ and is given by $q^2=G(p)$.} or, more explicitly, by a modular parameter $\tau\in\C$, $Im\tau>0$. Let us promote this modular parameter to an additional field $\tau(x)$ 
with Poisson brackets $\{\tau(x),\tau(y)\}=0,~~\{\tau(x),z_a(y)\}=2\pi i~ z_a(y)\delta^{\prime}(x-y)$. We will use the following

{\bf Lemma 1.} Consider a Poisson structure of the form:
\begin{equation}\label{homd}
 \{\tau(x),\tau(y)\}=0,~~\{\tau(x),z_a(y)\}=2\pi i~ z_a(y)\delta^{\prime}(x-y),
\end{equation}
$$\{z_a(x),z_b(y)\}=\sum_{1\leq c,d\leq n}q_{ab}^{cd}(\tau(x))z_c(x)z_d(x)\delta^{\prime}(x-y)+\sum_{1\leq c,d\leq n}r_{ab}^{cd}(\tau(x))z_c(x)z^{\prime}_d(x)\delta(x-y)+$$
$$+\sum_{1\leq c,d\leq n}t_{ab}^{cd}(\tau(x))z_c(x)z_d(x)\tau^{\prime}(x)\delta(x-y)$$
on the set of variables $\tau(x),~z_1(x),...,z_n(x)$ where $q_{ab}^{cd}(\tau),~r_{ab}^{cd}(\tau),~t_{ab}^{cd}(\tau)$ are functions in $\tau$. This Poisson structure can be descended to a family of Poisson structures on $\C P^{n-1}$ in the  following way:

{\bf 1.} Introduce affine coordinates on $\C P^{n-1}$ by $p_1(x)=\frac{z_1(x)}{z_n(x)},...,p_{n-1}(x)=\frac{z_{n-1}(x)}{z_n(x)}$.

{\bf 2.} Compute $\{p_a(x),p_b(y)\},~a,b=1,...,n-1$ by virtue of (\ref{homd}). Write these brackets as differential polynomials in $p_1(x),...,p_{n-1}(x)$.

{\bf 3.} Replace $\tau(x)$ by $\tau$ in obtained formulas for brackets $\{p_a(x),p_b(y)\}$ and declare that $\tau$ is a parameter.

{\bf Proof.} Note that $\{\tau(x),p_a(y)\}=\{\tau(x),\frac{z_a(x)}{z_n(x)}\}=0$ by virtue of (\ref{homd}). This means that the field $\tau(x)$ belongs to the center of the Poisson algebra generated by 
$\tau(x),~p_1(x),...,p_{n-1}(x)$ and can be set to a parameter which we also denote by $\tau$. It is also clear that we obtain a family of Poisson structures on $\C P^{n-1}$ (because we can replace $z_n(x)$ 
by any other $z_b(x)$ in our construction of affine coordinates $p_1(x),...,p_{n-1}(x)$). $\Box$

Using this Lemma we construct a lifting of the Poisson structure (\ref{phpr1}) as follows:

{\bf Proposition 2.} The following formulas
$$ \{\tau(x),\tau(y)\}=0,~~\{\tau(x),z_a(y)\}=2\pi i~ z_a(y)\delta^{\prime}(x-y),~a=1,2,$$

$$\{z_1(x),z_1(y)\}=\Big(-\frac{1}{6}g_2(\tau(x))z_1(x)z_2(x)+\frac{1}{2}g_3(\tau(x))z_2(x)^2\Big)\delta^{\prime}(x-y)+$$
$$+\delta(x-y)\Big(-\frac{1}{12}g_2(\tau(x))z_2(x)z_1^{\prime}(x)+(\frac{1}{2}g_3(\tau(x))z_2(x)-\frac{1}{12}g_2(\tau(x))z_1(x))z_2^{\prime}(x)+$$
$$+((6g_3(\tau(x))-4g_1(\tau(x))g_2(\tau(x)))z_1(x)z_2(x)+(18g_1(\tau(x))g_3(\tau(x))-g_2(\tau(x))^2)z_2(x)^2))\frac{i\tau^{\prime}(x)}{24\pi}\Big),$$ 

$$\{z_1(x),z_2(y)\}=\Big(2g_1(\tau(x))z_1(x)z_2(x)-\frac{1}{3}g_2(\tau(x))z_2(x)^2\Big)\delta^{\prime}(x-y)+$$
$$+\delta(x-y)\Big(g_1(\tau(x))z_2(x)z_1^{\prime}(x)+(g_1(\tau(x))z_1(x)-\frac{1}{3}g_2(\tau(x))z_2(x))z_2^{\prime}(x)-$$
$$-(2g_1(\tau(x))g_2(\tau(x))-3g_3(\tau(x)))z_2(x)^2\frac{i\tau^{\prime}(x)}{12\pi}\Big),$$ 

$$\{z_2(x),z_2(y)\}=\Big(-2z_1(x)z_2(x)+4g_1(\tau(x))z_2(x)^2\Big)\delta^{\prime}(x-y)+$$
$$+\Big(-z_2(x)z_1^{\prime}(x)+(4g_1(\tau(x))z_2(x)-z_1(x))z_2^{\prime}(x)+(12g_1(\tau(x))^2-g_2(\tau(x)))z_2(x)^2)\frac{i\tau^{\prime}(x)}{12\pi}\Big)\delta(x-y)$$ 
define a Poisson structure on the set of variables $\tau(x),~z_1(x),~z_2(x)$. Here the functions $g_1,~g_2,~g_3$ satisfy differential equations (\ref{dertau}). 

Moreover, the variables $p(x)=\frac{z_1(x)}{z_2(x)}$ satisfy (\ref{phpr1}) with $G(p)=-\frac{1}{2}(4p^3-g_2p-g_3)$. $\Box$

In this paper we construct a generalization of this homogeneous Poisson structure to an arbitrary number of variables $z_1(x),...,z_n(x)$. According to Lemma 1 this gives a family of Poisson structures of hydrodynamic type on $\C P^{n-1}$ for each $n=2,3,...$. This family for each $n$ is parametrized by a modular parameter $\tau$. 
 
In Section 2 we construct an analog of classical dynamical $r$-matrix which we call an $r$-matrix of hydrodynamic type. This is a certain Poisson structure of hydrodynamic type on the set of variables 
$\tau(x),~e(u,x)$ where $u$ is a spectral parameter and $\tau(x)$ plays the role of a dynamical variable. Structure constants of this Poisson structure can be written in terms of Jacobi forms with modular variable $\tau(x)$.

In Section 3 we explain how to construct homogeneous Poisson structures of hydrodynamic type on $\C^n$ for arbitrary $n$ using our $r$-matrix of hydrodynamic type. Our variable $e(u,x)$ serves as a generating function in this construction. Explicit (and really cumbersome) formulas for these homogeneous Poisson structures can be found in Appendix.

In Section 4 we collect definitions and results for elliptic functions and modular forms which we need in the main text.

In Section 5 we outline possible directions of further research.

\section{Elliptic dynamical \textit{r}-matrix of hydrodynamic type}
 
 Recall that a classical $r$-matrix $r_{ab}^{cd}(u,v)$ defines a Poisson structure of the forms\footnote{Here $u,~v$ are complex variables called spectral parameters. Notice that here and in the sequel  $r_{ab}^{cd}(u,v)$ are meromorphic functions in complex variables $u,~v$ which may have poles. Therefore, our Poisson structure is not defined at each values of spectral parameters. We require that anticommutativity and Jacobi identity hold each time when the corresponding Poisson brackets are defined.}
 
 \begin{equation}\label{r}
\{e_a(u),e_b(v)\}=\sum_{1\leq c,d\leq n}  r_{ab}^{cd}(u,v) e_c(u)e_d(v),~~~1\leq a,b\leq n.
 \end{equation}

 Moreover, a Poisson structure of the form (\ref{r}) on the set of variables $e_1(u),...,e_n(u)$  defines a classical $r$-matrix. One can also define a so-called dynamical $r$-matrix as a Poisson structure on the set 
 of variables $t_1,...,t_m,~e_1(u),...,e_n(u)$ of the form
 
 $$\{t_{\alpha},t_{\beta}\}=0,~~~\{t_{\alpha},e_a(u)\}=h_{\alpha a} e_a(u),$$
 \begin{equation}\label{rd}
\{e_a(u),e_b(v)\}=\sum_{1\leq c,d\leq n}  r_{ab}^{cd}(u,v,t_1,...,t_m)) e_c(u)e_d(v),~~~1\leq a,b\leq n,~1\leq \alpha,\beta\leq m.
 \end{equation}
 
We want to extend these definitions to the Poisson structures of hydrodynamic type and will call the corresponding object a (dynamical) $r$-matrix of hydrodynamic type.

Define brackets on the set of variables $\tau(x),~e(u,x)$ by\footnote{Here and in the sequel indexes like $e_u(u,x),~e_x(u,x)$ stand for partial derivatives.}

\begin{equation}\label{rh1}
\{\tau(x),\tau(y)\}=0,~~~\{\tau(x),e(u,y)\}=2\pi i~ e(u,y)~\delta^{\prime}(x-y), 
\end{equation}
\begin{equation}\label{rh2}
~~~~\{e(u,x),e(v,y)\}=
\end{equation}
$$=\Big(q(v,u,\tau(x))e_u(u,x)e(v,x)+q(u,v,\tau(x))e(u,x)e_v(v,x)+\lambda e_u(u,x)e_v(v,x)\Big)\delta^{\prime}(x-y)+$$
$$+\Big(q_{\tau}(v,u,\tau(x))e_u(u,x)e(v,x)\tau^{\prime}(x)+q_u(v,u,\tau(x))(e(u,x)e_x(v,x)-e(v,x)e_x(u,x))+$$
$$+q(v,u,\tau(x))e_u(u,x)e_x(v,x)+q(u,v,\tau(x))e(u,x)e_{vx}(v,x)+\lambda e_u(u,x)e_{vx}(v,x)\Big)\delta(x-y)$$
where\footnote{See (\ref{wz}), \ref{qg1}).} 
$$q(u,v,\tau)=\zeta(v-u,\tau)+\zeta(u,\tau)-g_1(\tau)v$$ 
and $\lambda$ is an arbitrary constant.

{\bf Theorem 1.} The brackets (\ref{rh1}), (\ref{rh2}) define a Poisson structure on the set of variables $\tau(x),~e(u,x)$.

{\bf Proof.} Using the identity (\ref{id}) we can write 
\begin{equation}\label{ql}
q(u,v,\tau)=\frac{1}{2}\cdot\frac{\wp_{u}(u,\tau)+\wp_{v}(v,\tau)}{\wp(v,\tau)-\wp(u,\tau)}+\zeta(v,\tau)-g_1(\tau)v 
\end{equation}
After substitution of this formula for $q$ into our brackets one can proof anticommutativity and Jacobi identity by direct calculation using formulas (\ref{derz}), (\ref{dertau}),  (\ref{derw})  for derivatives. $\Box$

The following statement provides infinitely many reductions of our Poisson structure (\ref{rh1}), (\ref{rh2}) in flat coordinates.

{\bf Theorem 2.} Define Poisson structure on the set of variables $t_1(x),...,t_{n-1}(x),\tau(x),f(x)$ as follows:
\begin{equation}\label{fl}
\{t_i(x),t_j(y)\}=\frac{1}{n}\delta^{\prime}(x-y),~~~\{t_i(x),t_i(y)\}=-\frac{n-1}{n}\delta^{\prime}(x-y),~~~1\leq i\ne j\leq n-1,
\end{equation}
$$\{\tau(x),f(y)\}=2\pi i~ f(y)~\delta^{\prime}(x-y),$$
$$\{\tau(x),\tau(y)\}=\{f(x),f(y)\}=\{\tau(x),t_i(y)\}=\{t_i(x),f(y)\}=0.$$
Then the variables\footnote{Here $\sigma(u,\tau)$ is the Weierstrass's sigma function, see (\ref{ws}),(\ref{ders1}), \ref{ders2}).}
$$e(u,x)=\frac{\sigma(u+\sum_{a=1}^{n-1}t_a(x),\tau(x))\prod_{a=1}^{n-1}\sigma(u-t_a(x),\tau(x))}{\sigma(u,\tau(x))^n}e^{-g_1(\tau(x))\sum_{1\leq a\leq b\leq n-1}t_a(x)t_b(x)}f(x)$$
satisfy (\ref{rh1}), (\ref{rh2}) with $\lambda=\frac{1}{n}$ by virtue of (\ref{fl}).  $\Box$

\section{Homogeneous Poisson structures on loop spaces}

Let $\mathcal F(z)$ be the space of elliptic functions in one variable $z$ with respect to $\Gamma=\{a+b\tau;~a,b\in\Z\}$ and holomorphic outside $\Gamma$.
  For $n\in\N$ let $\mathcal F_n(z)\subset\mathcal F(z)$ be the subspace of functions with poles of order $\leqslant n$ on $\Gamma$.
  It is clear that $\dim\mathcal F_n(z)=n$. It is known that the functions 
  $$p_{2\alpha}(z,\tau)=\wp(z,\tau)^{\alpha},~~~ p_{2\alpha+3}(z,\tau)=
  -\frac{1}{2}\wp(z,\tau)^{\alpha}\wp_z(z,\tau),~~~ \alpha=0,1,2,...$$ 
  form a basis of the vector space $\mathcal F(z)$ and the functions $p_0(z,\tau),~p_2(z,\tau),p_3(z,\tau),...,p_n(z,\tau)$ form a basis of the space $\mathcal F_n(z)$. Note that $p_{\alpha}(z,\tau)$ has a pole of order $\alpha$ at $z=0$.
 
 {\bf Theorem 3.} Set $\lambda=\frac{1}{n}$ and
 \begin{equation}\label{red}
 e(u,x)=p_0(u,\tau(x))z_0(x)+\sum_{2\leq a\leq n}p_a(u,\tau(x))z_a(x) 
 \end{equation}
 in (\ref{rh1}), (\ref{rh2}). Then the formulas (\ref{rh1}), (\ref{rh2}) define a homogeneous Poisson structure of hydrodynamic type on the set of variables $\tau(x),~z_0(x),z_2(x),...,z_n(x)$.
 
 {\bf Proof.} Substitution of (\ref{red}) into (\ref{rh1}) gives the following brackets:
 \begin{equation}\label{hom1}
 \{\tau(x),\tau(y)\}=0,~~~\{\tau(x),z_a(y)\}=2\pi i~ z_a(y)~\delta^{\prime}(x-y),~~~a=0,2,3,...,n. 
 \end{equation}
Let us substitute (\ref{red}) into (\ref{rh2}). In the l.h.s. we obtain
$$\sum_{a,b=0,2,...,n}\Big(p_a(u,\tau(x))p_b(v,\tau(y))\{z_a(x),z_b(y)\}+2\pi i(p_a(u,\tau(x))p_b(v,\tau(y)))_{\tau}z_a(x)z_b(y)\delta^{\prime}(x-y)\Big)$$
Let us denote by $\Delta(u,v)$ the difference of the l.h.s. and the r.h.s. of (\ref{rh2}). One can check by direct calculation that $\Delta(u,v)$ is an elliptic function with respect to each of variable $u,~v$. Moreover, $\Delta(u,v)$ belongs to $\mathcal F_n(u)$ as a function of $u$ and it belongs to $\mathcal F_n(v)$ as a function of $v$.  Therefore, $\Delta(u,v)$ is an element of the vector space $\mathcal F_n(u)\otimes \mathcal F_n(v)$. It is clear that $\{p_a(u,\tau)p_b(v,\tau);~a,b=0,2,...,n\}$ is a basis\footnote{Using the property  $\delta(x-y)F(x)=\delta(x-y)F(y)$ where $F$ is an arbitrary function we replace all terms containing $\tau(y)$ by the corresponding terms with $\tau(x)$. After that we identify $\tau(x)\equiv\tau$.} of $\mathcal F_n(u)\otimes \mathcal F_n(v)$. Expanding $\Delta(u,v)$ by this basis and equating to zero coefficient at $p_a(u,\tau)p_b(v,\tau)$ we obtain a formula for $\{z_a(x),z_b(y)\}$. $\Box$
 
{\bf Remark 1.} Explicit formulas for $\{z_a(x),z_b(y)\}$ are cumbersome, see Appendix. To obtain these formulas we write $e(u,x)$ in the form
$$e(u,x)=e_0(\wp(u,\tau(x)),x)-\frac{1}{2}\wp_u(u,\tau(x))e_1(\wp(u,\tau(x)),x).$$
Comparing with (\ref{red}) we see that $e_0(z,x),~e_1(z,x)$ are polynomials with respect to $z$. Let us also use the formula (\ref{ql}) for $q(u,v,\tau)$. Substituting these expressions for $e$ and $q$ into 
the difference of the l.h.s. and the r.h.s. of (\ref{rh2}) and using formulas (\ref{derz}), (\ref{dertau}),  (\ref{derw}) to eliminate derivatives by $\tau$ and to reduce powers of $\wp_u(u,\tau),~\wp_v(v,\tau)$ to either 0 or 1 we obtain a polynomial in $\wp_u(u,\tau),~\wp_v(v,\tau)$ of degree 1 with respect to each variable. Equating to zero coefficients of this polynomial we obtain explicit formulas from Appendix.

 \section{Weierstrass's functions and modular forms}
 
 Here we collect formulas which we need in the main text. Proofs (in slightly different notations) can be found in \cite{lang,zagier}. Let $\tau\in\C$ be a modular parameter, we assume $Im(\tau)>0$. 
 
 The Weierstrass's elliptic function is defined by
 
 \begin{equation}\label{wp}
 \wp(z,\tau)=\frac{1}{z^2}+\sum_{(a,b)\in\Z^2\setminus \{0\}}\Big(\frac{1}{(z-a-b\tau)^2}-\frac{1}{(a+b\tau)^2}\Big).
 \end{equation}
 
 The Weierstrass zeta function is defined by
 
 \begin{equation}\label{wz}
 \zeta(z,\tau)=\frac{1}{z}+\sum_{(a,b)\in\Z^2\setminus \{0\}}\Big(\frac{1}{z-a-b\tau}+\frac{1}{a+b\tau}+
  \frac{z}{(a+b\tau)^2}\Big). 
 \end{equation}
 
 The following formulas for derivatives of these functions with respect to $z$ hold:
 
 \begin{equation}\label{derz}
 \zeta_z(z,\tau)=-\wp(z,\tau),~~ \wp_z(z,\tau)^2=4\wp(z,\tau)^3-g_2(\tau)\wp(z,\tau)-g_3(\tau).
 \end{equation}
 
 Here the modular forms $g_2(\tau),~g_3(\tau)$ are given by
 
 \begin{equation}\label{qg23}
g_2(\tau)=\frac{4}{3}\pi^4+320\pi^4\sum_{m=1}^{\infty} \frac{m^3q^m}{1-q^m},~~~g_3(\tau)=\frac{8}{27}\pi^6-\frac{448}{3}\pi^6\sum_{m=1}^{\infty} \frac{m^5q^m}{1-q^m},~~~q=e^{2\pi i\tau}.
\end{equation}

 Notice that the Weierstrass's zeta function $\zeta(z,\tau)$ is not elliptic. In particular, we have
 
 \begin{equation}\label{g1}
 \zeta(z+1,\tau)=\zeta(z,\tau)+g_1(\tau). 
 \end{equation}

where $g_1(\tau)$ is given by

\begin{equation} \label{qg1}
g_1(\tau)=\frac{1}{3}\pi^2-8\pi^2\sum_{m=1}^{\infty} \frac{mq^m}{1-q^m},~~~\\ q=e^{2\pi i\tau}.
\end{equation}

The following formulas for derivatives with respect to $\tau$ hold:

$$2\pi ig_1^{\prime}(\tau)=\frac{1}{12}g_2(\tau)-g_1(\tau)^2,$$
\begin{equation}\label{dertau}
2\pi ig_2^{\prime}(\tau)=6g_3(\tau)-4g_1(\tau)g_2(\tau),
\end{equation}

$$2\pi ig_3^{\prime}(\tau)=\frac{1}{3}g_2(\tau)^2-6g_1(\tau)g_3(\tau).$$

Notice that we need the function $g_1(\tau)$ in order to write these formulas for derivatives. Moreover, derivatives of Weierstrass's functions can also be computed:

\begin{equation}\label{derw}
2\pi i\zeta_{\tau}(z,\tau)=g_1(\tau)z\wp(z,\tau)+\frac{1}{12}g_2(\tau)z-g_1(\tau)\zeta(z,\tau)-\zeta(z,\tau)\wp(z,\tau)-\frac{1}{2}\wp_z(z,\tau), 
\end{equation}

$$2\pi i\wp_{\tau}(z,\tau)=(\zeta(z,\tau)-zg_1(\tau))\wp_z(z,\tau)+2\wp(z,\tau)^2-2g_1(\tau)\wp(z,\tau)-\frac{1}{3}g_2(\tau).$$
 
We also need the following identity:

\begin{equation}\label{id}
\zeta(z_1-z_2,\tau)-\zeta(z_1,\tau)+\zeta(z_2,\tau)=\frac{1}{2}\cdot\frac{\wp_{z_1}(z_1,\tau)+\wp_{z_2}(z_2,\tau)}{\wp(z_1,\tau)-\wp(z_2,\tau)}.
\end{equation}

Proof of the identities (\ref{derw}), (\ref{id}) is standard: to check that the difference between the r.h.s. and the l.h.s is an elliptic function and to calculate the decomposition at each pole.

The Weierstrass's sigma function is defined by
 
 \begin{equation}\label{ws}
 \sigma(z,\tau)=z\prod_{(a,b)\in\Z^2\setminus \{0\}}\Big(1-\frac{z}{a+b\tau}\Big)e^{\frac{z}{a+b\tau}+\frac{z^2}{2(a+b\tau)^2}}.
 \end{equation}
 
The following formulas for derivatives with respect to $z$ and $\tau$ hold: 

\begin{equation}\label{ders1}
\frac{\sigma_z(z,\tau)}{\sigma(z,\tau)}=\zeta(z,\tau), 
\end{equation}

\begin{equation}\label{ders2}
2\pi i \frac{\sigma_{\tau}(z,\tau)}{\sigma(z,\tau)}=g_1(\tau)+\frac{1}{24}g_2(\tau)z^2-zg_1(\tau)\zeta(z,\tau)+\frac{1}{2}\zeta(z,\tau)^2-\frac{1}{2}\wp(z,\tau).
\end{equation}

\section{Conclusion and outlook}

It follows from (\ref{red}) that the space of generators of our homogeneous Poisson structure on $\C^n$ with a basis $z_0,z_2,...,z_n$ can be identified with a dual space to a space of theta functions 
of order $n$. Therefore, our Poisson structure is invariant with respect to a projective action of $(\Z/n\Z)^2$. It will be interesting to write our Poisson structure in a different basis, so that it become manifestly $(\Z/n\Z)^2-$invariant. More generally, it will be interesting to investigate if there exist other $(\Z/n\Z)^2-$invariant homogeneous Poisson structures of hydrodynamic type. For example, in the case of usual homogeneous Poisson structures there exist Poisson algebras $q_{n,k}(\tau)$ parametrized by a modular parameter $\tau$ and a discrete parameter $k\in(\Z/n\Z)^*$ \cite{fod,od1,od2}.

Our $r$-matrix of hydrodynamic type is connected with a moduli space of elliptic curves. It will be interesting to generalize this to higher genus.

It will be interesting to investigate a quantization (if any) of our Poisson structures of hydrodynamic type, for example, a quantization of our $r$-matrix.

\section{Appendix}

The brackets $\{z_a(x),z_b(y)\}$ can be written as follows:
$$\{z_a(x),z_b(y)\}=P_{ab}~\delta^{\prime}(x-y)+Q_{ab}~\delta(x-y)$$
where $P_{ab}$ are homogeneous quadratic polynomials in $z_0(x),z_2(x),...,z_n(x)$ and $Q_{ab}$ are homogeneous quadratic polynomials in $z_0(x),z_2(x),...,z_n(x),~z^{\prime}_0(x),z^{\prime}_2(x),...,z^{\prime}_n(x)$. Let 
$$e_0(u,x)=\sum_{0\leq a\leq \frac{n}{2}}u^az_{2a}(x),~~~e_1(u,x)=\sum_{0\leq a\leq \frac{n-3}{2}}u^az_{2a+3}(x)$$
We have

$$\sum_{0\leq a,b\leq \frac{n}{2}}P_{2a,2b}u^av^b=2ue_{0u}(u)e_0(v)g_1-\frac{1}{6}(-12u^2v+g_2u+2g_2v+3g_3)\frac{e_{0u}(u)e_0(v)}{u-v}+2ve_{0v}(v)e_0(u)g_1+$$
$$+\frac{1}{6}(-12uv^2+2g_2u+g_2v+3g_3)\frac{e_{0v}(v)e_0(u)}{u-v}+\frac{1}{8}(4v^3-g_2v-g_3)(4u^3-g_2u-g_3)\frac{e_{1u}(u)e_1(v)}{u-v}-$$
$$-\frac{1}{8}(4v^3-g_2v-g_3)(4u^3-g_2u-g_3)\frac{e_{1v}(v)e_1(u)}{u-v}+\frac{1}{4n}(4v^3-g_2v-g_3)(4u^3-g_2u-g_3)e_{1u}(u)e_{1v}(v)+$$
$$-(3u^2v^2-\frac{1}{4}g_2(u-v)^2+\frac{1}{16}g_2^2+\frac{3}{4}g_3u+\frac{3}{4}g_3v)e_1(u)e_1(v)+\frac{1}{8n}(12v^2-g_2)(4u^3-g_2u-g_3)e_{1u}(u)e_1(v)+$$
$$+\frac{1}{8n}(12u^2-g_2)(4v^3-g_2v-g_3)e_{1v}(v)e_1(u)+\frac{1}{16n}(12v^2-g_2)(12u^2-g_2)e_1(u)e_1(v)$$

$$\sum_{\substack{0\leq a\leq \frac{n}{2}\\0\leq b\leq \frac{n-3}{2}}}P_{2a,2b+3}u^av^b=2ue_{0u}(u)e_1(v)g_1+\frac{1}{6}(12u^2v-g_2u-2g_2v-3g_3)\frac{e_{0u}(u)e_1(v)}{u-v}+$$
$$+\frac{1}{2}(4u^3-g_2u-g_3)\frac{e_{1u}(u)e_0(v)-e_{0v}(v)e_1(u)}{u-v}+2ve_{1v}(v)e_0(u)g_1+3e_0(u)e_1(v)g_1-$$
$$-\frac{1}{6}(12uv^2-2g_2u-g_2v-3g_3)\frac{e_{1v}(v)e_0(u)}{u-v}-\frac{1}{4}(12uv-g_2)\frac{e_0(u)e_1(v)}{u-v}+\frac{1}{4}(12u^2-g_2)\frac{e_0(v)e_1(u)}{u-v}+$$
$$+\frac{1}{n}(4u^3-g_2u-g_3)e_{0v}(v)e_{1u}(u)+\frac{1}{n}(6u^2-\frac{1}{2}g_2)e_{0v}(v)e_1(u)$$

$$\sum_{0\leq a,b\leq \frac{n-3}{2}}P_{2a+3,2b+3}u^av^b=2ue_{1u}(u)e_1(v)g_1+6e_1(u)e_1(v)g_1+\frac{4}{n}e_{0u}(u)e_{0v}(v)+$$
$$+\frac{1}{6}(12u^2v-g_2u-2g_2v-3g_3)\frac{e_{1u}(u)e_1(v)}{u-v}-\frac{1}{6}(12uv^2-2g_2u-g_2v-3g_3)\frac{e_{1v}(v)e_1(u)}{u-v}+$$
$$+2\frac{e_0(v)e_{0u}(u)-e_0(u)e_{0v}(v)}{u-v}+2ve_{1v}(v)e_1(u)g_1$$

$$\sum_{0\leq a,b\leq \frac{n}{2}}Q_{2a,2b}u^av^b=\frac{1}{6}(12g_1uv-12g_1v^2-12uv^2+2g_2u+g_2v+3g_3)\frac{e_0(u)e_{0vx}(v)}{u-v}-$$
$$-(4v^3-g_2v-g_3)(4u^3-g_2u-g_3)\frac{e_1(u)e_{1vx}(v)}{8(u-v)}+(12g_2g_1u-g_2^2+18g_1g_3-18g_3u)\frac{ie_0(u)e_{0v}(v)\tau^{\prime}(x)}{24\pi(u-v)}+$$
$$+(12g_1u^2-12g_1uv+12u^2v-g_2u-2g_2v-3g_3)\frac{e_{0u}(u)e_{0x}(v)}{6(u-v)} +$$
$$+\frac{e_0(v)e_{0x}(u)-e_0(u)e_{0x}(v)}{4(u-v)^2}(4g_1(u-v)^2+4u^2v+4uv^2-g_2u-g_2v-2g_3)-$$
$$-\frac{ie_1(v)e_{1u}(u)\tau^{\prime}(x)}{8\pi n}(2g_2g_1-3g_3)(4u^3-g_2u-g_3)+\frac{1}{8}(4v^3-g_2v-g_3)(4u^3-g_2u-g_3)\frac{e_{1u}(u)e_{1x}(v)}{u-v}-$$
$$-\frac{ie_1(v)e_{1u}(u)\tau^{\prime}(x)}{48\pi(u-v)}(4u^3-g_2u-g_3)(12g_2g_1v-g_2^2+18g_1g_3-18g_3v)+$$
$$+\frac{1}{8}(4v^3-g_2v-g_3)(4u^3-g_2u-g_3)\frac{e_1(v)e_{1x}(u)}{(u-v)^2}-\frac{ie_1(u)e_1(v)\tau^{\prime}(x)}{16\pi n}(2g_2g_1-3g_3)(12u^2-g_2)-$$
$$-\frac{e_1(u)e_{1x}(v)}{16(u-v)^2}(48u^4v^2-64u^3v^3+48u^2v^4-4g_2u^4+8g_2u^3v-24g_2u^2v^2+8g_2uv^3-4g_2v^4+g_2^2u^2+g_2^2v^2+$$
$$+4g_3u^3-12g_3u^2v-12g_3uv^2+4g_3v^3+2g_2g_3u+2g_2g_3v+2g_3^2)+$$
$$+\frac{ie_0(v)e_{0u}(u)\tau^{\prime}(x)}{24\pi(u-v)}(24g_1^2u^2-24g_1^2uv-8g_2g_1u-4g_2g_1v-2g_2u^2+2g_2uv+g_2^2-18g_1g_3+12g_3u+6g_3v)+$$
$$+\frac{e_{1u}(u)e_{1vx}(v)}{4n}(4v^3-g_2v-g_3)(4u^3-g_2u-g_3)+\frac{e_1(u)e_{1vx}(v)}{8n}(12u^2-g_2)(4v^3-g_2v-g_3)-$$
$$-\frac{ie_{1u}(u)e_{1v}(v)\tau^{\prime}(x)}{24\pi n}(4u^3-g_2u-g_3)(12g_2g_1v-g_2^2+18g_1g_3-18g_3v)+$$
$$+\frac{ie_1(u)e_{1v}(v)\tau^{\prime}(x)}{48\pi(u-v)}(4u^3-g_2u-g_3)(12g_2g_1v-g_2^2+18g_1g_3-18g_3v)-$$
$$-\frac{ie_1(u)e_{1v}(v)\tau^{\prime}(x)}{48\pi n}(12u^2-g_2)(12g_2g_1v-g_2^2+18g_1g_3-18g_3v)+\frac{e_1(u)e_{1x}(v)}{16n}(12v^2-g_2)(12u^2-g_2)+$$
$$+\frac{e_{1u}(u)e_{1x}(v)}{8n}(12v^2-g_2)(4u^3-g_2u-g_3)+\frac{ie_1(u)e_1(v)\tau^{\prime}(x)}{48\pi}(24g_2g_1u^2-24g_2g_1uv-6g_1g_2^2+$$
$$+4g_2^2u+2g_2^2v-72g_1g_3u-36g_1g_3v-36g_3u^2+36g_3uv+9g_2g_3) $$

$$\sum_{\substack{0\leq a\leq \frac{n}{2}\\0\leq b\leq \frac{n-3}{2}}}Q_{2a,2b+3}u^av^b=-\frac{1}{2}(4u^3-g_2u-g_3)\frac{e_1(u)e_{0vx}(v)}{u-v}+$$
$$+\frac{1}{6}(12g_1uv-12g_1v^2-12uv^2+2g_2u+g_2v+3g_3)\frac{e_0(u)e_{1vx}(v)}{u-v}+$$
$$+\frac{1}{6}(12g_1u^2-12g_1uv+12u^2v-g_2u-2g_2v-3g_3)\frac{e_{0u}(u)e_{1x}(v)}{u-v}+$$
$$+\frac{ie_1(v)e_{0u}(u)\tau^{\prime}(x)}{24\pi} (24g_1^2u^2-24g_1^2uv-8g_2g_1u-4g_2g_1v-2g_2u^2+2g_2uv+g_2^2-18g_1g_3+12g_3u+6g_3v)+$$
$$+\frac{1}{4(u-v)^2}e_1(v)e_{0x}(u)(4g_1u^2-8g_1uv+4g_1v^2+4u^2v+4uv^2-g_2u-g_2v-2g_3)+$$
$$+\frac{1}{2}(4u^3-g_2u-g_3)\frac{e_{0x}(v)e_{1u}(u)}{u-v}+\frac{1}{4}(4u^3-12u^2v+g_2u+g_2v+2g_3)\frac{e_{0x}(v)e_1(u)}{(u-v)^2}+$$
$$+\frac{1}{2}(4u^3-g_2u-g_3)\frac{e_0(v)e_{1x}(u)}{(u-v)^2}+\frac{ie_0(u)e_{1v}(v)\tau^{\prime}(x)}{24\pi(u-v)}(12g_2g_1u-g_2^2+18g_1g_3-18g_3u)+$$
$$+\frac{1}{2}(4g_1u^2-8g_1uv+4g_1v^2-8u^2v+4uv^2+g_2u+g_3)\frac{e_{1x}(v)e_0(u)}{(u-v)^2}+$$
$$+\frac{i(e_0(u)e_1(v)-e_1(u)e_0(v))\tau^{\prime}(x)}{24\pi(u-v)^2}(12g_1g_2u-g_2^2+18g_1g_3-18g_3u)+$$
$$+\frac{1}{n}(4u^3-g_2u-g_3)e_{0vx}(v)e_{1u}(u)+\frac{1}{2n}(12u^2-g_2)e_{0vx}(v)e_1(u)$$

$$\sum_{0\leq a,b\leq \frac{n-3}{2}}Q_{2a+3,2b+3}u^av^b=\frac{1}{6}(12g_1uv-12g_1v^2-12uv^2+2g_2u+g_2v+3g_3)\frac{e_1(u)e_{1vx}(v)}{u-v}+$$
$$+2\frac{e_0(v)e_{0x}(u)-e_0(u)e_{0x}(v)}{(u-v)^2}+\frac{1}{6}(12g_1u^2-12g_1uv+12u^2v-g_2u-2g_2v-3g_3)\frac{e_{1u}(u)e_{1x}(v)}{u-v}+$$
$$+\frac{ie_1(v)e_{1u}(u)\tau^{\prime}(x)}{24\pi(u-v)} (24g_1^2u^2-24g_1^2uv-8g_2g_1u-4g_2g_1v-2g_2u^2+2g_2uv+g_2^2-18g_1g_3+12g_3u+6g_3v)+$$
$$+\frac{1}{4}(4g_1(u-v)^2+4u^2v+4uv^2-g_2u-g_2v-2g_3)\frac{e_1(v)e_{1x}(u)}{(u-v)^2}+2\frac{e_{0u}(u)e_{0x}(v)-e_0(u)e_{0vx}(v)}{u-v}+$$
$$+\frac{ie_1(u)e_{1v}(v)\tau^{\prime}(x)}{24\pi(u-v)}(12g_2g_1u-g_2^2+18g_1g_3-18g_3u)+\frac{ie_1(u)e_1(v)\tau^{\prime}(x)}{8\pi}(12g_1^2-g_2)+$$
$$+\frac{1}{4}(20g_1(u-v)^2-4u^2v-4uv^2+g_2u+g_2v+2g_3)\frac{e_{1x}(v)e_1(u)}{(u-v)^2} +\frac{4}{n}e_{0vx}(v)e_{0u}(u)$$

Here the r.h.s. do not depend on $y$ and we suppress dependence of $x$, for example
$$\tau\equiv \tau(x),~~~e_0(u)\equiv e_0(u,x),~~~e_1(u)\equiv e_1(u,x).$$
We also suppress dependence of $\tau$ in $g_1,~g_2,~g_3$. Recall that indexes $u,~v,~x$ stand for partial derivatives.

\end{document}